\documentclass[12pt, a4paper, twoside]{amsart}

\usepackage[T1]{fontenc}
\usepackage[applemac]{inputenc}
\usepackage{amsmath}
\usepackage{amsthm}
\usepackage{amsfonts}
\usepackage{bbm}
\usepackage[english]{babel}
\usepackage{enumerate}
\usepackage{bm}
\usepackage{ae}
\usepackage{amssymb}

\theoremstyle{definition}

\theoremstyle{remark}

\theoremstyle{plain}
\newtheorem{thm}{Theorem}
\newtheorem{lem}[thm]{Lemma}

\newtheorem{cor}[thm]{Corollary}

\makeatletter
\newcommand*{\house}[1]{%
  \mathord{%
    \mathpalette\@house{#1}%
  }%
}
\newcommand*{\@house}[2]{%
  \dimen@=\fontdimen8 %
      \ifx#1\scriptscriptstyle\scriptscriptfont
      \else\ifx#1\scriptstyle\scriptfont
      \else\textfont\fi\fi
      3 %
  \sbox0{%
    $#1%
      \vrule width\dimen@\relax
      \overline{%
        \kern2\dimen@
        \begingroup % to keep changes of \dimen@ in #2 local
          #2%
        \endgroup
        \kern2\dimen@
      }%
      \vrule width\dimen@\relax
      \mathsurround=1.5\dimen@ % outside margin
    $%
  }%
  \ht0=\dimexpr\ht0-\dimen@\relax
  \dp0=\dimexpr\dp0+2\dimen@\relax
  \vbox{%
    \kern\dimen@ % reinsert previously removed space
    \copy0 %
  }%
}
\makeatother

\begin{document}

\title[Continued fractions]{Irrationality and transcendence of continued fractions with algebraic integers}

\author{SIMON B. ANDERSEN}

\address{S. B. Andersen, Aarhus Private Gymnasium, Bjødstrupvej 26, DK-8270 Højbjerg, Denmark}

\email{simon\_andersen94@hotmail.com}

\author{SIMON KRISTENSEN}

\address{S. Kristensen, Department of Mathematics, Aarhus University, Ny Munkegade 118,
  DK-8000 Aarhus C, Denmark}

\email{sik@math.au.dk}

\begin{abstract}
We extend a result of Han\v{c}l, Kolouch and Nair on the irrationality and transcendence of continued fractions. We show that for a sequence $\{\alpha_n\}$ of algebraic integers of bounded degree, each attaining the maximum absolute value among their conjugates and satisfying certain growth conditions, the condition 
$$
\limsup_{n \rightarrow \infty} \vert\alpha_n\vert^{\frac{1}{Dd^{n-1} \prod_{i=1}^{n-2}(Dd^i + 1)}} = \infty
$$
implies that the continued fraction $\alpha = [0;\alpha_1, \alpha_2, \dots]$ is not an algebraic number of degree less than or equal to $D$.
\end{abstract}


\maketitle

\section{Introduction}

Continued fractions are a convenient tool in studying the arithmetical properties of a number. Usually, one considers a real number $\alpha$ and the simple continued fraction
$$
\alpha = a_0 + \cfrac{1}{a_1 + \cfrac{1}{a_2 + \cfrac{1}{\ddots}}} = [a_0; a_1, a_2, \dots],
$$
where $a_0 \in \mathbb{Z}$ and $a_i \in \mathbb{N}$ for $i \ge 1$. The sequence of partial quotients $\{a_i\}$ may be finite or infinite. A number is rational if and only if the sequence is finite. 

For an irrational number $\alpha$, we define the sequence of convergents of $\alpha$ to be the finite continued fractions obtained by truncation,
$$
\frac{p_n}{q_n} = [a_0; a_1, a_2, \dots, a_n].
$$
It is a well-known consequence of Roth's theorem \cite{MR1451873} that $\alpha$ is transcendental if
$$
\limsup_{n \rightarrow \infty} \frac{\log a_{n+1}}{\log q_n} > 0.
$$

While choosing the partial quotients to be natural numbers is entirely natural for real numbers, this is in fact not necessary. For instance, in Hurwitz' continued fractions \cite{MR1554754} the partial quotients are Gaussian integers, and any complex number can be expanded in a Hurwitz continued fraction. In fact, one may use as partial quotients any choice of non-zero real or complex numbers, though the resulting continued fraction need not be convergent. If the partial quotients are all positive real numbers, the convergence is ensured by the condition 
$$
\sum_{n=1}^\infty a_n = \infty,
$$
see \cite{MR1451873} where it is shown that this is necessary and sufficient in this case. For more general complex partial quotients this is no longer the case, but if the sequence $\{\vert q_n\vert\}$ is increasing, this is sufficient for the convergence of the continued fraction.

In the present note, we are interested in continued fractions where the partial quotients are algebraic integers. In particular, we are interested in general criteria for the irrationality or transcendence of the resulting continued fraction. Not much appears to be known in this setting. In \cite{MR3826638}, a growth condition on a series of positive integers $\{a_n\}$ is given, which guarantees that the continued fraction 
$$
\gamma = [0;\sqrt{a_1}, \sqrt{a_2}, \dots]
$$
is not algebraic of degree $D$. In brief, it is shown that if
$$
\limsup_{n \rightarrow \infty} a_n^{\frac{1}{D2^{n-1} \prod_{i=1}^{n-2}(D2^i + 1)}} = \infty,
$$
then $\gamma$ cannot be an algebraic number of degree at most $D$. An immediate corollary is that the condition
$$
\limsup_{n \rightarrow \infty} a_n^{2^{-n^2}} = \infty
$$
implies the transcendence of $\gamma$.

In the present note -- much in the spirit of \cite{Andersen:2018aa} -- we extend this result to deal with sequences of algebraic integers. To state the result, recall first that the \emph{house} of an algebraic number $\alpha$ is the maximal modulus among the conjugates of $\alpha$, i.e. if the conjugates of $\alpha$ are $\alpha_1 = \alpha, \alpha_2, \dots, \alpha_d$, then the house of $\alpha$ is defined to be
$$
\house{\alpha} = \max_{1 \le j \le d} \{\vert \alpha_j\vert\}.
$$

Next, recall that for a sequence of non-zero complex numbers $\alpha_n$, the denominators $q_n$ of the convergents of the continued fraction $[0;\alpha_1, \alpha_2, \dots]$ satisfy the recurrence
\begin{equation}
\label{eq:CF}
q_0 = 1, \quad q_1 = \alpha_1, \quad q_{n+2}= \alpha_{n+2} q_{n+1} + q_n,
\end{equation}
see \cite{MR1451873}.

\begin{thm}
\label{thm:main}
Let $d, D \in \mathbb{N}$, such that $d>1$, and let $\{\alpha_n\}$ be a sequence of algebraic integers with $\max \deg \alpha_n = d$, such that $\vert \alpha_n \vert = \house{\alpha_n}$, and such that
$$
\limsup_{n \rightarrow \infty} \house{\alpha_n}^{\frac{1}{Dd^{n-1} \prod_{i=1}^{n-2}(Dd^i + 1)}} = \infty.
$$
Let $\alpha$ be given by the continued fraction  $\alpha=[0;\alpha_1,\alpha_2,...]$, and suppose that the sequence of absolute values $\{\vert q_n\vert\}$, where $q_n$ is defined by \eqref{eq:CF}, is increasing. Then,
$$
\deg\left(\alpha\right) > D.
$$
\end{thm}

The final assumption on the recurrence sequence is only to ensure that the $q_n$ are well defined and that the continued fraction converges. If the $\alpha_n$ are all real, positive numbers, we can completely remove the condition on the recurrence sequence, as it is trivially satisfied.

\begin{cor}
\label{cor:main}
Let $d, D \in \mathbb{N}$, such that $d>1$, and let $\{\alpha_n\}$ be a sequence of real and positive algebraic integers with $\max \deg \alpha_n = d$, such that $\vert \alpha_n \vert = \house{\alpha_n}$, and such that
$$
\limsup_{n \rightarrow \infty} \house{\alpha_n}^{\frac{1}{Dd^{n-1} \prod_{i=1}^{n-2}(Dd^i + 1)}} = \infty.
$$
Let $\alpha$ be given by the continued fraction $\alpha=[0;\alpha_1,\alpha_2,...]$. Then,
$$
\deg\left(\alpha\right) > D.
$$
\end{cor}

Finally, we remark that if
$$
\limsup_{n \rightarrow \infty}  \house{\alpha_n}^{d^{-n^2}} = \infty,
$$
the growth conditions in both Theorem \ref{thm:main} and Corollary \ref{cor:main} are automatically satisfied for any value of $D \in \mathbb{N}$. Hence, we obtain the following immediate corollary.

\begin{cor}
Let $\{\alpha_n\}$ be a sequence of algebraic integers satisfying the conditions of either Theorem \ref{thm:main} or Corollary \ref{cor:main}. If in addition 
$$
\limsup_{n \rightarrow \infty}  \house{\alpha_n}^{d^{-n^2}} = \infty,
$$
the continued fraction $\alpha=[0;\alpha_1,\alpha_2,...]$ is transcendental.
\end{cor}

\section{Auxiliary results}

For the proof of Theorem \ref{thm:main}, we will need several auxillary results, which will be stated without proof. Their proofs can be found in the references given. The first is classical and relates the Weil height $H(\alpha)$ of an algebraic number $\alpha$ to it's Mahler measure $M(\alpha)$, see e.g. \cite{MR1756786}.

\begin{thm}
\label{thm:height_measure}
For an algebraic number $\alpha$ of degree $d$,
$$
H(\alpha) = M(\alpha)^{1/d}.
$$
\end{thm}

In \cite{Andersen:2018aa}, we use this to derive the following lemma.

\begin{lem}
\label{lem:height_house_measure}
Let $\alpha$ be an algebraic integer of degree $d$. Then,
$$
H(\alpha) = M(\alpha)^{1/d} \le \house{\alpha} \le M(\alpha) = H(\alpha)^d.
$$
The inequalities are best possible.
\end{lem}

That the Weil height of an algebraic number is invariant under taking the reciprocal is also classical, see e.g. \cite{MR1756786}. We phrase this as a lemma.

\begin{lem}
\label{lem:reciprocal}
Let $\alpha$ be a non-zero algebraic number. Then, $H(\alpha) = H(1/\alpha)$.
\end{lem}

We will need bounds on the height and degree of sums of algebraic numbers. The first inequality in the next lemma is proved in \cite{MR1756786}, and the second inequality can be found in  \cite{MR0258803}.

\begin{lem}
\label{lem:elementary}
Let $n \in \mathbb{N}$, and let $\beta_1, \dots, \beta_n$ be algebraic numbers. Then,
$$
H\left(\sum_{i=1}^n \beta_i\right) \le 2^n \prod_{i=1}^n H(\beta_i),
$$
and
$$
\deg\left(\sum_{i=1}^n \beta_i\right) \le \prod_{i=1}^n \deg(\beta_i).
$$
\end{lem}

We will also need the Liouville--Mignotte bounds \cite{liouville, MR554376} on the distance between algebraic numbers. A unified proof of the inequalities can be found in  \cite[Appendix A]{MR2136100}

\begin{lem}
\label{lem:liouville}
Let $\alpha$ and $\beta$ be non-conjugate algebraic numbers. Then,
$$
\vert \alpha - \beta \vert \ge \frac{1}{2^{\deg(\alpha)\deg(\beta)}M(\alpha)^{\deg(\beta)}M(\beta)^{\deg(\alpha)}}.
$$
\end{lem}

Again from \cite{Andersen:2018aa}, we will use the following result.

\begin{lem}
\label{lem:calculus2}
Let $\{a_n\}_{n=1}^\infty$ be a sequence of real numbers such that 
$$
\limsup_{n \rightarrow \infty} a_n = \infty.
$$
Then for infinitely many $k \in \mathbb{N}$,
$$
a_{k+1} > \left(1+\frac{1}{k^2}\right) \max_{1 \le n \le k} a_n.
$$
\end{lem}

Finally, from the theory of continued fractions we will need the following result, see \cite{MR1451873}

\begin{lem}
\label{lem:continued fraction}
Let $\alpha_n\geq 1$ for all $n\in\mathbb{N}$. For the continued fraction $\alpha=[0;\alpha_1,\alpha_2,...]$ set $\alpha(n)=[0;\alpha_1,...,\alpha_{n}]$. Then
$$\left|\alpha-\alpha(n)\right|<\frac{1}{|q_{n+1}q_n|},$$
where $q_0=1$, $q_1=\alpha_1$, and $q_n=\alpha_nq_{n-1}+q_{n-2}$.
\end{lem}

\section{Proof of the main theorem}

First, Corollary \ref{cor:main} immediately follows from Theorem \ref{thm:main}, as Lemma \ref{lem:height_house_measure} implies that $\alpha_n = \house{\alpha_n} \ge 1$, since $M(\alpha_n) \ge 1$ trivially. Consequently, the sequence $\{q_n\}$ is evidently increasing.

We now give a  proof of Theorem \ref{thm:main}. We argue by contradiction, so suppose that $\alpha$ is algebraic with $\deg(\alpha) \le D$. For $N \in \mathbb{N}$, we let 
$$
\alpha(N) = \frac{p_N}{q_N} = [0;\alpha_1, \alpha_2, \dots, \alpha_N].
$$
Note that this is an algebraic number, and that 
$$
[0;\alpha_1, \alpha_2, \dots, \alpha_N]^{-1} = \alpha_1 + \frac{1}{[0;\alpha_2, \alpha_3, \dots, \alpha_N]}.
$$
We estimate the degree and the Mahler measure of $\alpha(N)$ using this and Lemma \ref{lem:elementary}. For the degree, we get
\begin{multline*}
\deg(\alpha(N)) = \deg([0;\alpha_1,...,\alpha_N]) \\
= \deg([0;\alpha_1,...,\alpha_N]^{-1}) \le  d\deg([0;\alpha_2,...,\alpha_N])
\le \dots\le d^N
\end{multline*}
For the Mahler measure, we get 
\begin{alignat*}{2}
M(\alpha(N)) &\le H(\alpha(N))^{d^N} \\
&=H(\alpha(N)^{-1})^{d^N} \\
&\le (2^2H(\alpha_1)H([0;\alpha_2,...,\alpha_N]))^{d^N} \\
&\le ...\le \left(2^{2(N-1)} \prod_{n=1}^N H \left(\alpha_n\right)\right)^{d^N} \le \left(2^{2(N-1)} \prod_{n=1}^N  \house{\alpha_n}\right)^{d^N},
\end{alignat*}
where the latter equality follows from Lemma \ref{lem:height_house_measure}. 

Now, using the above estimates and Lemma \ref{lem:liouville},
\begin{alignat*}{2}
\vert \alpha - \alpha(N) \vert &\ge \frac{1}{2^{\deg(\alpha)\deg(\alpha(N))}M(\alpha)^{\deg(\alpha(N))}M(\alpha(N))^{\deg(\alpha)}}\\
&\ge \frac{1}{2^{Dd^N}M(\alpha)^{d^N} \left(2^{2(N-1)} \prod _{n=1}^N  \house{\alpha_n}\right)^{Dd^N}}\\
&= \frac{1}{\left(2^{2N -1}H(\alpha) \prod _{n=1}^N  \house{\alpha_n}\right)^{Dd^N}}.
\end{alignat*}
The upshot is the following critical estimate, valid for all $N \in \mathbb{N}$,
\begin{equation}
\label{eq:critical}
\vert\alpha-\alpha(N)\vert \left(2^{2N -1}H(\alpha) \prod _{n=1}^N  \house{\alpha_n}\right)^{Dd^N}\ge 1.
\end{equation}

By Lemma \ref{lem:calculus2}, for infinitely many values of $N$, 
$$
 \house{\alpha_{N+1}}^{\frac{1}{Dd^N\prod_{i=1}^{N-1}(Dd^i + 1)}} \ge \left(1 + \frac{1}{N^2}\right) \max_{1 \le n \le N}  \house{\alpha_n}^{\frac{1}{Dd^{n-1}\prod_{i=1}^{n-2}(Dd^i + 1)}}.
 $$
 Note that for these infinitely many values of $N$, we get from Lemma \ref{lem:continued fraction} that for $N$ large enough,
 \begin{equation}\label{eq:inequality}
 \vert \alpha-\alpha(N)\vert<\frac{1}{\vert q_{N+1}q_N \vert}\le \frac{1}{\house{\alpha_{N+1}}}.
 \end{equation}
 The final inequality follows by the recursion \eqref{eq:CF} and the assumption that $\house{\alpha_{N+1}} = \vert \alpha_{N+1}\vert$, which in turn implies that $\vert \alpha_{N+1}\vert \ge 1$ by Lemma \ref{lem:height_house_measure}.
 Note also that
 $$
 \log\left(1 + \frac{1}{N^2}\right) \ge \frac{2N^2 - 1}{2N^4},
 $$
 and that
 $$
 d^N\left(\prod_{n=1}^{N-1}(Dd^n + 1)\right) \frac{2N^2 - 1}{2N^4} \ge \log(2) (d+1)^N,
 $$
 so that 
 $$
 \left(1 + \frac{1}{N^2}\right)^{Dd^N\prod_{i=1}^{N-1}(Dd^i+1)} > 2^{D(d+1)^N}.
 $$
 
 Hence, we find that for these infinitely many values of $N$,
 \begin{alignat*}{2}
 \house{\alpha_{N+1}} &\ge 2^{D(d+1)^N} \\
 &\quad \left(\max_{1 \le n \le N}  \house{\alpha_n}^{\frac{1}{Dd^{n-1}\prod_{i=1}^{n-2}(Dd^i + 1)}}\right)^{Dd^N\prod_{i=1}^{N-1}(Dd^i +1)} \\
 & =  2^{D(d+1)^N}  \\
 &\quad \left(\max_{1 \le n \le N}  \house{\alpha_n}^{\frac{1}{Dd^{n-1}\prod_{i=1}^{n-2}(Dd^i +1)}}\right)^{D^2d^{2N-1}\prod_{i=1}^{N-2}(Dd^i +1)} \\
  &\quad \left(\max_{1 \le n \le N}  \house{\alpha_n}^{\frac{1}{Dd^{n-1}\prod_{i=1}^{n-2}(Dd^i + 1)}}\right)^{Dd^N\prod_{i=1}^{N-2}(Dd^i +1)} \\
  & \ge  2^{D(d+1)^N}\house{\alpha_N}^{Dd^N}  \\
  &\quad \left(\max_{1 \le n \le N}  \house{\alpha_n}^{\frac{1}{Dd^{n-1}\prod_{i=1}^{n-2}(Dd^i + 1)}}\right)^{Dd^N\prod_{i=1}^{N-2}(Dd^i +1)} \\
 &\ge \cdots \ge \\
 &\ge 2^{D(d+1)^N}\left(\prod_{n=1}^N \house{\alpha_n}\right)^{D d^N}.
 \end{alignat*}
Inserting this together with \eqref{eq:inequality} into \eqref{eq:critical}, we find that
 $$
 1 \le \vert\alpha-\alpha(N)\vert \left(2^{2N -1}H(\alpha) \prod _{n=1}^N  \house{\alpha_n}\right)^{Dd^N}
$$
$$
\le \dfrac{\left(2^{2N -1}H(\alpha) \prod _{n=1}^N  \house{\alpha_n}\right)^{Dd^N}}{\house{\alpha_{N+1}}}\le\dfrac{\left(2^{2N -1}H(\alpha) \prod _{n=1}^N  \house{\alpha_n}\right)^{Dd^N}}{2^{D(d+1)^N}\left(\prod_{n=1}^N \house{\alpha_n}\right)^{D d^N}}
$$
$$
=\dfrac{\left(2^{2N -1}H(\alpha) \right)^{Dd^N}}{2^{D(d+1)^N}},
$$
which is not true for large enough $N$.
\qed

\providecommand{\bysame}{\leavevmode\hbox to3em{\hrulefill}\thinspace}
\providecommand{\MR}{\relax\ifhmode\unskip\space\fi MR }
\providecommand{\MRhref}[2]{%
  \href{http://www.ams.org/mathscinet-getitem?mr=#1}{#2}
}
\providecommand{\href}[2]{#2}

\end{document}